\documentclass[11pt,twoside]{article}
\usepackage{amsthm}
\usepackage{amsmath}
\usepackage{amssymb}
\usepackage{array}
\usepackage{graphicx}
\usepackage{graphics}
\usepackage{epsfig}
\usepackage{amsfonts,amscd}
\usepackage{exscale}
\usepackage{eucal}
\usepackage{latexsym,amsmath}
\usepackage{indentfirst}
\numberwithin{equation}{section}
\newtheorem{Prop}{\bf Proposition}[section]

\newtheorem{Rem}{\bf Remark}[section]
\newtheorem{Ex}{\bf Example}[section]
\newtheorem{Th}{Theorem}[section]

\theoremstyle{definition} \theoremstyle{plain}

\begin{document}
\def \b{\Box}

\begin{center}
{\Large {\bf Matrix algorithm for determination of the elementary\\[0.2cm]
paths and elementary circuits using exotic semirings}}
\end{center}

\begin{center}
{\bf Gheorghe IVAN}
\end{center}

\setcounter{page}{1}

\pagestyle{myheadings}

{\small {\bf Abstract}. We propose a new method for determining
the elementary paths and elementary circuits in a directed graph.
Also, the Hamiltonian paths and Hamiltonian circuits are
enumerated.}
{}

\section{Introduction}
\smallskip
\indent\indent Idempotent mathematics is based on replacing the
usual arithmetic operations with a new set of basic operations,
that is on replacing numerical fields by idempotent semirings.
Exotic semirings such as the max-plus algebra ${\bf R}_{max}$ or
concatenation semiring ${\cal P}(\Sigma^{\ast})$ have been
introduced in connection with various fields: graph theory, Markov
decision processes, language theory, discrete event systems
theory, see \cite{kusa}, \cite{fink}, \cite{mohr}, \cite{bagu}.

In this paper we will have to consider various semirings, and will
universally use the notation $\oplus, \otimes, \varepsilon, e$
with a context dependent manning (e.g. $\oplus:= max $ in ${\bf
R}_{max}$ but $\oplus:= \cup$ in ${\cal P}(\Sigma^{\ast})$,
$\varepsilon:= - \infty $ in ${\bf R}_{max}$ but $\varepsilon:=
\emptyset$ in ${\cal P}(\Sigma^{\ast})$).

 In many fields of applications, the graphs are used
widely for modelling of practical problems. This paper will focus
on two algebraic path problems, namely: the {\it elementary path
problem} ({\bf EPP}) and the {\it elementary circuit problem}
({\bf ECP}).

For a directed graph $G=(V,E)$ ($|V|=n$) and $u,v\in V$, ({\bf
EPP}) and ({\bf ECP}) are formulated as follows:\\[0.1cm]
$\bullet~~({\bf EPP})$ {\it enumerate the elementary paths from
$u$ to $v$ of length $k$ ($1\leq k\leq n-1$)};\\[0.1cm]
$\bullet~~({\bf ECP})$ {\it enumerate the elementary circuits
starting in $u$ of length $k$ ($1\leq k\leq n$)}.

For to solve the above problems we give a method based on a
$n\times n$ matrix with entries in the semiring of distinguished
languages.

 These algebraic path problems are applied into large
domains: combinatorial optimization, traffic control, Internet
routing etc.

The paper is organized as follows. In Section 2 we construct a
special idempotent semiring denoted by ${\cal
P}^{\ast}(\Sigma_{dw}^{\ast})$ and named the semiring of
distinguished languages. The semiring of matrices with entries in
${\cal P}^{\ast}(\Sigma_{dw}^{\ast})$ is presented in Section 3.
This algebraic tool is used to establishing of a one-to-one
correspondence between the distinguished words and the elementary
paths in a directed graph. In Section 4 we give an algorithm to
determine all the elementary paths and elementary circuits in a
directed graph. This new practical algorithm is based on the latin
composition of distinguished languages. Another method for
determination of the elementary paths is the well-known latin
multiplication technique of Kaufmann (see \cite{kauf}).

\section{Semiring of distinguished formal languages}
\indent\indent We start this section by recalling of some
necessary backgrounds on semirings for our purposes (see
\cite{mohr}, \cite{bagu}, \cite{litv} and references therein for
more details).

{\bf Semirings}.  Let $S$ be a nonempty set endowed with two
binary operations, {\it addition} (denoted with $\oplus$) and {\it
multiplication} (denoted with $\otimes$). The algebraic structure
$(S,\oplus, \otimes, \varepsilon, e )$ is a \textit{semiring}, if
it fulfills the following conditions:

$(1)~~(S,\oplus, \varepsilon)~$ is a commutative monoid with
$\varepsilon$ as the neutral element for $\oplus;$

$(2)~(S, \otimes,  e )~$ is a  monoid with $\varepsilon$ as the
identity element for $\otimes;$

$(3)~~\otimes $ distributes over $\oplus;$

$(4)~~\varepsilon~$ is an absorbing element for $\otimes$, that is
$~a\otimes \varepsilon=\varepsilon \otimes a= \varepsilon,~\forall
a\in S.$

A semiring where addition is idempotent (that is, $~a\oplus
a=a,~\forall a \in S$) is called an {\it idempotent semiring}. If
$\otimes$ is commutative, we say that $S$ is a {\it commutative
semiring}.

In the following we introduce the {\it monoid of distinguished
words over an alphabet}.

An {\it alphabet} is a finite set $\Sigma$ of symbols. An {\it
word} over $\Sigma$ is a finite sequence of symbols of the
alphabet $\Sigma$. The number of symbols of a word $x$ is called
the {\it length of $x$} and its length is represented by $ |x|$.
The {\it empty word}, denoted with $\lambda$, is the word with
length zero (that is, it has no symbols).

A {\it simple word of length $k$} with $1\leq k\leq n$ over the
alphabet $\Sigma$ ($|\Sigma|=n$) is a word of the form $a = a_{1}
a_{2}...a_{k-1}a_{k}$ such that $a_{i}\neq a_{j}$ for $i\neq j $
and $i,j=\overline{1,k}$.

A {\it simple cyclic word of length $k+1$} with $1\leq k\leq n$
over $\Sigma$ is a word of the form $b = a_{1} a_{2}\ldots
a_{k-1}a_{k}a_{1}$ such that $a = a_{1} a_{2}\ldots a_{k-1}a_{k}$
is a simple word. The simple cyclic words of length $2$ over
$\Sigma$ are the words of the form $u = a a$ for $a \in \Sigma$.

By a {\it distinguished word} over $\Sigma$ we mean a word $w$
over $\Sigma$ which is a simple word of length $k$ or a simple
cyclic word of length $k+1$ with $1\leq k\leq n$. Denote the {\it
set of distinguished words over $\Sigma$ completed with the word}
$\lambda$  by $\Sigma_{dw}^{\ast}$.

\markboth{Gheorghe Ivan}{Matrix algorithm for determination of the
elementary paths and . . .}

For example, if $\Sigma = \{a,b\}$, then
$~\Sigma_{dw}^{\ast}=\{\lambda, a, b, aa, bb, ab, ba, aba, bab\}.$

It is easy to prove that, {\it if $~\Sigma$ is an alphabet with
$|\Sigma|=n$ symbols, then $\Sigma_{dw}^{\ast}$ is a finite set
with $\sigma_{n}$ elements, where}\\[-0.2cm]
\begin{equation}
\sigma_{n} = 1 + 2 n! + 2 \sum\limits_{k=1}^{n-1}\label{(2.1)}
\frac{n!}{(n-k)!}.
\end{equation}

On the set $\Sigma_{dw}^{\ast}$ we introduce the binary operation
$\circ_{\ell}$ given as follows:\\[0.1cm]
$(1)~$ for all distinguished word $ x \in \Sigma_{dw}^{\ast}$ and
simple cyclic word $c\in \Sigma_{dw}^{\ast}$, we have
 \begin{equation}
\lambda \circ_{\ell} x = x\circ_{\ell}\lambda =
\lambda~~~\hbox{and}~~~ c \circ_{\ell} x = x\circ_{\ell} c =
\lambda;\label{(2.2)}
\end{equation}
$(2)~~$ Let  $x = a_{1} a_{2}\ldots a_{k-1} a_{k}$ and $y = b_{1}
b_{2}\ldots b_{r-1} b_{r}$  be two simple words of the lengths $k$
and $r$ with $1\leq k, r\leq n$. The word $x \circ_{\ell} y \in
\Sigma_{dw}^{\ast} $ is defined by:\\[-0.1cm]
\begin{equation}
x \circ_{\ell} y = \left \{
\begin{array}{ll}
a_{1}\ldots a_{k}b_{2}\ldots b_{r-1} b_{r},&
\hbox{if}~~a_{k}= b_{1}, \{a_{1},\ldots,  a_{k}\}\cap\{b_{2},\ldots, b_{r}\}= \emptyset\\
a_{1}\ldots a_{k}b_{2}\ldots b_{r-1} a_{1}, & \hbox{if}~~ a_{k}=
b_{1}, b_{r}=a_{1}, \{a_{1},\ldots,
a_{k}\}\cap\{b_{2},\ldots, b_{r-1}\}= \emptyset\\
\lambda,& \hbox{otherwise}.\label{(2.3)}
\end{array}\right.
\end{equation}

The operation $\circ_{\ell}$ is called the {\it latin composition
of distinguished words}.
\begin{Ex}
{\rm The set $\Sigma_{dw}^{\ast}$ of distinguished words over
$\Sigma=\{1, 2, 3, 4\}$ has $\sigma_{4}=129$ elements. If $ x =
123, y=31, z=1, c_{1}= 22$ and
$c_{2}= 343$, then \\[0.1cm]
$c_{1}\circ_{l} x = 22 \circ_{l} 123 =\lambda,~~~x \circ_{l} c_{2}
= 123 \circ_{l} 343 =\lambda,~~~ x \circ_{l} y = 123 \circ_{l} 31
= 1231,$\\[0.1cm]
$y \circ_{l} x = 31 \circ_{l} 123 =3123,~~~ z \circ_{l} x = 1
\circ_{l} 123 = 123.~$ Note that $ x \circ_{l} y \neq y \circ_{l}
x$}.\hfill$\Box$
\end{Ex}

For all $a, b, c\in \Sigma_{dw}^{\ast}$ we have $~( a \circ_{\ell}
b ) \circ_{\ell} c = a \circ_{\ell} ( b \circ_{\ell} c )$. Then $(
\Sigma_{dw}^{\ast}, \circ_{\ell}, \lambda )$ is a monoid, called
the {\it monoid of distinguished words generated by alphabet
$\Sigma$}.

A {\it distinguished language} $L$ over alphabet $\Sigma$ is a
subset of distinguished words over $\Sigma$  completed with the
empty word $\lambda$, that is $L\subseteq \Sigma_{dw}^{\ast}$ and
$\lambda \in L$.

{\bf Convention}. $(i)~$ If $L=\{\lambda\}$ or $L=\{a\}$ where $a$
is a distinguished word, then we shall use the notations $\lambda$
and $a$, respectively.

$(ii)~$ If $L\neq \{\lambda \}$, then we enumerate only its
distinguished words of length $k\geq 1. \hfill\Box$

The monoid $\Sigma_{dw}^{\ast}$ contains two special distinguished
languages, namely: the language $\{\lambda\}$ (it contains only
the empty word $\lambda$) and the language $\Sigma$ (it contains
all symbols of the alphabet and the empty word $\lambda$).

The {\it set of distinguished languages over an alphabet $\Sigma$}
is ${\cal P}^{\ast}(\Sigma_{dw}^{\ast}).$

Since distinguished languages are sets, all the set operations can
be applied to distinguished languages. Then the union and
intersection of two distinguished languages are distinguished
languages.

The {\it latin composition} of the distinguished languages $L_{1}$
and $L_{2}$ is the distinguished language $L_{1} \circ_{\ell}
L_{2}$ defined by
\begin{equation}
L_{1} \circ_{\ell} L_{2} =\{ x\circ_{\ell} y~|~ x \in L_{1}
~\hbox{and}~ y\in L_{2} \},\label{(2.4)}
\end{equation}
that is, $L_{1} \circ_{\ell} L_{2}$ is the set of distinguished
words obtaining by the latin composition of words of $L_{1}$ with
those of $L_{2}$.

For example, if $L_{1},L_{2}\subseteq \Sigma_{dw}^{\ast}$ are
distinguished languages over
$\Sigma=\{1, 2, 3, 4\}$, where $L_{1}= \{2, 412\}$ and $L_{2}= \{11, 23\}$, then\\[0.1cm]
$L_{1}\circ_{\ell} L_{2}= \{2, 412\}\circ_{\ell}\{11, 23\} =
 \{ 2\circ_{\ell} 11, 412 \circ_{\ell} 11, 2 \circ_{\ell} 23,
412 \circ_{\ell} 23\}= \{ 23, 4123\}$.

\begin{Prop}
$(i)~$ Let $ L, L_{1}, L_{2}, L_{3}\in {\cal
P}(\Sigma_{dw}^{\ast}).$ Then:
\begin{equation}
L \circ_{l}\lambda = \lambda \circ_{l} L=\lambda,~~~  L
\circ_{l}\Sigma = \Sigma \circ_{l} L= L,~~~ ( L_{1} \circ_{l}
L_{2}) \circ_{l} L_{3}=  L_{1} \circ_{l} ( L_{2} \circ_{l}
L_{3}).\label{(2.5)}
\end{equation}

$(ii)~$ The set ${\cal P}^{\ast}(\Sigma_{dw}^{\ast}) $ endowed
with multiplication $\circ_{l}$ has a structure of monoid.
\end{Prop}
\smallskip
{\it Proof.} $(i)~$ Using the definitions one easily verify that
$(2.5)$ holds.

$(ii)$ From $(i)$ it follows that $({\cal
P}^{\ast}(\Sigma_{dw}^{\ast}), \circ_{l}, \Sigma ) $ is a monoid.
\hfill$\Box$\\

On the set ${\cal P}^{\ast}(\Sigma_{dw}^{\ast}) $ of distinguished
languages we define the binary operations:
\begin{equation}
L_{1} \oplus L_{2}:= L_{1} \cup L_{2}~~~\hbox{and}~~~L_{1} \otimes
L_{2}:= L_{1} \circ_{\ell} L_{2},~~~\forall L_{1}, L_{2} \in {\cal
P}^{\ast}(\Sigma_{dw}^{\ast}). \label{(2.6)}
\end{equation}

\begin{Prop}
$({\cal P}^{\ast}(\Sigma_{dw}^{\ast}), \cup, \circ_{\ell},
\lambda, \Sigma )$ is an idempotent semiring.
\end{Prop}
{\it Proof.}  For to verify the conditions from definition of an
idempotent semiring we apply the properties of the union of sets
and Proposition 2.1. \hfill$\Box$

We call $({\cal P}^{\ast}(\Sigma_{dw}^{\ast}),\cup, \circ_{\ell},
\lambda, \Sigma ) $ the {\it semiring of distinguished languages
over $\Sigma$}.

\section{Matrices over semirings and directed graphs}

Let $(S, \oplus, \otimes, \varepsilon, e )$ be an (idempotent)
semiring.  For each positive integer $n$, let $ M_{n}(S)$ be
denote the set of $n\times n$ matrices with entries in $S$. The
operations $\oplus$ and $\otimes$ on $S$ induce corresponding
operations on $ M_{n}(S)$ in the obvious way. Indeed, if
$A=(A_{ij}), B=(B_{ij})   \in M_{n}(S)$ then we have:\\[-0.2cm]
\[
A\oplus B= ((A\oplus B)_{ij}) ~~~\hbox{and}~~~ A\otimes B=
((A\otimes B)_{ij}),~~~i,j =\overline{1,n}~~~\hbox{where}
\]\\[-0.8cm]
\begin{equation}
(A\oplus B)_{ij}:= A_{ij} \oplus B_{ij}~~~\hbox{and}~~~ (A\otimes
B)_{ij}:= \bigoplus\limits_{k=1}^{n} A_{ik}\otimes
B_{kj}.\label{(3.1)}
\end{equation}

The set $M_{n}(S)$ contains two special matrices with entries in
$S$, namely the zero matrix $O_{\oplus n}$, which  has all its
entries equal to $\varepsilon$, and the identity matrix
$I_{\otimes n}$, which  has the diagonal entries equal to $e$ and
the other entries equal to $\varepsilon$.

It is easy to check that the following proposition holds.
\begin{Prop}
$( M_{n}(S), \oplus, \otimes, O_{\oplus n}, I_{\otimes n}) $ is an
idempotent semiring, where the operations $\oplus$ and $\otimes$
are given in $(3.1)$.\hfill$\Box$
\end{Prop}

We call $( M_{n}(S), \oplus, \otimes, O_{\oplus n}, I_{\otimes n})
$ the {\it semiring of $n\times n$ matrices with entries in $S$}.
In particular, if $S:= ({\cal P}^{\ast}(\Sigma_{dw}^{\ast}),\cup,
\circ_{\ell}, \emptyset, \Sigma )$, then  $( M_{n}({\cal
P}^{\ast}(\Sigma_{dw}^{\ast})), \cup, \circ_{\ell}, O_{\oplus n},
I_{\otimes n}) $ is called the {\it semiring of $n\times n$
matrices over $ {\cal P}^{\ast}(\Sigma_{dw}^{\ast})$}. The
operation $\otimes:=\circ_{\ell} $ is called the {\it
multiplication of matrices based on latin composition of words}.
\begin{Ex}
{\rm  Let be semiring  $( M_{2}({\cal
P}^{\ast}(\Sigma_{dw}^{\ast})), \cup, \circ_{\ell}, O_{\oplus 2},
I_{\otimes 2}) $ with $\Sigma = \{a, b, c\}.$ The product
$A\circ_{\ell} B$ of the $A, B $ with entries in the semiring
${\cal P}^{\ast}(\Sigma_{dw}^{\ast})$ is}
\[
A\circ_{\ell} B= \left(\begin{array}{cc} ab & \varepsilon\\
bca & bc\\
\end{array}\right)\circ_{\ell}\left(\begin{array}{cc}
b &ab\\
c & \varepsilon\\
\end{array}\right) =\left(\begin{array}{cc}
(ab \circ_{\ell} b)\oplus (\varepsilon \circ_{\ell} c) & (ab
\circ_{\ell} ab)\oplus (\varepsilon \circ_{\ell} \varepsilon) \\
(bca \circ_{\ell} b)\oplus (bc \circ_{\ell} c) & (bca
\circ_{\ell} ab)\oplus (bc \circ_{\ell} \varepsilon) \\
\end{array}\right)=
\]
\[
=\left(\begin{array}{cc} ab \oplus \varepsilon & \varepsilon
\oplus \varepsilon  \\
\varepsilon\oplus bc & bcab \oplus \varepsilon) \\
\end{array}\right )= \left(\begin{array}{cc} ab  \cup\varepsilon  & \varepsilon
\cup \varepsilon  \\
\varepsilon \cup bc & bcab \cup  \varepsilon) \\
\end{array}\right )= \left(\begin{array}{cc} ab & \varepsilon
\\
bc & bcab  \\
\end{array}\right)
.~~~\hfill\Box
\]
\end{Ex}
 A {\it directed graph} is a pair $ G=(V,E)$ where $V$ is a finite
{\it set of vertices} of the graph $G$  and $E\subseteq V\times V$
is a {\it set of arcs} of $G$. A typical arc $(u,v)\in E$ is
thought of as an arrow directed from $u$ to $v$.

Let $G=(V,E)$ be a directed graph with $|V|=n$. A {\it path from
$u$ to $v$ of length $k$} ($k\geq 1$) in $G$ is a sequence of
vertices $p =( v_{1}, v_{2},\ldots, v_{k}, v_{k+1})$ with $
v_{1}=u, v_{k+1}=v$ such that $(v_{i}, v_{i+1})\in E$ for all
$i=\overline{1,k}$; $v_{1}$ is called the {\it starting vertex}
 and $v_{k+1}$ the {\it end-vertex} of $p$,
respectively. The length of path $p$ will denoted by $\ell(p)$.

A path  $p =( v_{1}, v_{2},\ldots, v_{k}, v_{k+1})$ is called {\it
circuit} if $ v_{k+1}= v_{1}$ and $k\geq 1$. In particular, for
$k=1$ we obtain the circuit $(v_{1}, v_{1})$ of length $1$.

 We denote with $P(v_{i}, v_{j}, k)$ ($k\geq 1$),  the  set of all  paths of length $k$ from the
 starting vertex $v_{i}\in V$ to end-vertex $v_{j}\in V$. In particular, when
 $v_{i}=v_{j}$, $C(v_{i}, k)=P(v_{i}, v_{i}, k)$ (
$k\geq 1$) is the set of all circuits of length $k$ starting at
vertex $v_{i}$.

A path  $p =( v_{1}, v_{2},\ldots, v_{k}, v_{k+1})$  is called an
{\it elementary path}  from $v_{1}$ to $v_{k+1}$, if $k\geq 1$ and
$ v_{i}\neq v_{j}$ for $i\neq j$ and $i,j=\overline{1,k+1}.$ A
circuit $c =( u_{1}, u_{2},\ldots, u_{k}, u_{1})$ with $\ell(c)=k$
is called an {\it elementary circuit}, if $( u_{1}, u_{2},\ldots,
u_{k})$ is an elementary path.

 We denote with $P_{elem}(v_{i}, v_{j}, k)$ ($k\geq 1$),  the  set of all elementary paths
 of length $k$ from  $v_{i}\in V$ to  $v_{j}\in V$. In particular, when $v_{i}=v_{j}$, then
 $C_{elem}(v_{i}, k)=P_{elem}(v_{i}, v_{i}, k)$
($k\geq 1$) is the set of all elementary circuits of length $k$
starting at $v_{i}$.

A {\it Hamiltonian path} (resp., {\it circuit}) is a path (resp.,
circuit) that contains each vertex exactly once. Hence, a
Hamiltonian path (resp., circuit) is an elementary path $p_{H}$
with $\ell(p_{H})=n-1$ (resp., an elementary circuit $c_{H}$ with
$\ell(c_{H})=n$).

A {\it weighted directed graph} is a graph $G=(V,E)$ with a
mapping $w: E\rightarrow S$ that assigns each arc $(u,v)\in E$ a
weight $w(u,v)$ from the semiring $(S, \oplus, \otimes,
\varepsilon, e)$. A weighted directed graph with the cost function
$w$ is denoted by $G=(V,E,w)$.

The {\it weight} or {\it cost} of path $p =( v_{1},\ldots, v_{k},
v_{k+1})$ is the element $w(p)\in S$ where
\begin{equation}
w(p)=\bigotimes\limits_{i=1}^{k} w(v_{i}, v_{i+1}). \label{(3.1)}
\end{equation}

To each given weighted directed graph $G=(V, E, w)$ with
$V=\{v_{1}, v_{2}, \ldots, v_{n}\}$ we can associate a  $n\times n
$ matrix $M_{w}(G)$ with entries in a semiring $(S, \oplus,
\otimes , \varepsilon , e)$ as follows. For this, we define the
matrix $M_{w}(G)= (M_{ij})\in M(n, S)$ where
\begin{equation}
M_{ij}=\left \{\begin{array}{lll} w(v_{i}, v_{j}) & \hbox{if} &
(v_{i}, v_{j})\in E\\
\varepsilon & \hbox{if} &
(v_{i}, v_{j})\notin E\\
\end{array}\right. \label{(3.2)}
\end{equation}

{\it To each directed graph $G=(V, E)$ with $V=\{v_{1}, v_{2},
\ldots, v_{n}\}$ we can associate two weight functions} in the
following way.

$\bullet ~~$ Let be the {\it numbering semiring} $({\bf N}, +,
\cdot, 0, 1)$ of natural numbers endowed with the usual addition
and multiplication. Consider the weight function $w_{a}: E\to {\bf
N}$ defined by $ w_{a}(v_{i}, v_{j})=1 $ for all $(v_{i},
v_{j})\in E.$ The matrix $M_{w_{a}}(G)\in M(n, {\bf N})$, denoted
with $A$, is called the {\it adjacency matrix} of graph $G$.

$\bullet ~~$ Let be the idempotent semiring $({\cal
P}^{\ast}(\Sigma_{dw}^{\ast}), \cup, \circ_{\ell}, \emptyset,
\Sigma )$ of distinguished languages over alphabet $\Sigma =V$.
Define the weight function $w_{\ell}: E\to {\cal
P}^{\ast}(\Sigma_{dw}^{\ast})$ given by $ w_{\ell}(v_{i}, v_{j})=
v_{i}v_{j} $ for all $(v_{i}, v_{j})\in E$ (that is,
$w_{\ell}(v_{i}, v_{j})$  is the distinguished language which
contains only the distinguished word $v_{i}v_{j}$ of length $1$).
The matrix $M_{w_{\ell}}(G)\in M(n, {\cal
P}^{\ast}(\Sigma_{dw}^{\ast}))$ is denoted with $L$ and is called
the {\it latin matrix} of $G$.

More precisely, the adjacency matrix $A=(A_{ij})\in M(n, {\bf N})$
and  latin matrix $L=(L_{ij})\in M(n, {\cal
P}^{\ast}(\Sigma_{dw}^{\ast}))$ associated to graph $G$ are given
by
\begin{equation}
A_{ij}=\left \{\begin{array}{lll} 1 & \hbox{if} &
(v_{i}, v_{j})\in E\\
0 & \hbox{if} &
(v_{i}, v_{j})\notin E\\
\end{array}\right.~~~\hbox{and}~~~L_{ij}=\left \{\begin{array}{lll} v_{i}v_{j} & \hbox{if} &
(v_{i}, v_{j})\in E\\
\varepsilon & \hbox{if} &
(v_{i}, v_{j})\notin E\\
\end{array}\right. \label{(3.3)}
\end{equation}

\section{Matrix algorithm for enumerating of elementary paths
and elementary circuits in a directed  graph}

Consider the latin matrix $L=(L_{ij})\in M(n, {\cal
P}^{\ast}(\Sigma_{dw}^{\ast}))$ associated to graph $G=(V,E)$
($|V|=n$), defined by $(3.4)$. Using the multiplication of
matrices based on latin composition, we define by recurrence the
power $L^{[k]}$ of matrix $L$ in the following
way:\\[-0.2cm]
\begin{equation}
L^{[1]}=L,~~ L^{[2]}=
L\circ_{\ell}L,~~\ldots,~~L^{[k]}=L\circ_{\ell}L^{[k-1]}~~\hbox{for}~~k\geq
2 .\label{(4.1)}
\end{equation}

Applying $(3.1)$ and replacing  $\oplus$ and $\otimes$ with the
correspondent operations of the semiring $ M(n, {\cal
P}^{\ast}(\Sigma_{dw}^{\ast}))$, we have $L^{[k]}=(L_{ij}^{[k]})$
for $i,j=\overline{1,n}$, where\\[-0.2cm]
\begin{equation}
L_{ij}^{[k]}=\bigcup\limits_{m=1}^{n}(L_{im}\circ_{\ell}
L_{mj}^{[k-1]}),~~k\geq 2.\label{(4.2)}
\end{equation}

It is easy to prove that:\\[-0.2cm]
\begin{equation}
L^{[n]} ~~\hbox{ is a diagonal matrix}~~\hbox{and}~~
L^{[n+q]}=O_{\oplus n}~~\hbox{for all}~~q\geq 1.\label{(4.3)}
\end{equation}
\begin{Th}
Let $G=(V,E)$ be a directed graph with  $V=\{v_{1}, v_{2}, \ldots,
v_{n}\}$ and the latin matrix $L=(L_{ij})\in M(n, {\cal
P}^{\ast}(\Sigma_{dw}^{\ast}))$. If $L^{[k]}=(L_{ij}^{[k]})$,
then:\\[-0.2cm]
\begin{equation}
L_{ij}^{[k]} = P_{elem}(v_{i}, v_{j}, k),~~~~
i,j=\overline{1,n},~i \neq j,~1\leq k\leq n-1; \label{(4.4)}
\end{equation}
\begin{equation}
L_{ii}^{[k]} = C_{elem}(v_{i}, k),~~~~~~~~~~~~~~~~~~~~~~
i=\overline{1,n},~1\leq k\leq n. \label{(4.5)}
\end{equation}
\end{Th}

{\it Proof.} We proceed by induction on $k$. Choose $v_{i}$ and
$v_{j}$ arbitrarily.

{\bf (1)~Case $i\neq j$}.  For $k=1$, the relation $(4.4)$ holds.
Indeed, we have that the only elementary path in $P_{elem}(v_{i},
v_{j}, 1)$ is the arc $(v_{i}, v_{j})$, since $L_{ij}^{[1]}=
v_{i}v_{j}.$ Therefore, $P_{elem}(v_{i}, v_{j}, 1)= \{ (v_{i},
v_{j})\}$. Note that if $L_{ij}=\varepsilon,$ then
$P_{elem}(v_{i}, v_{j}, 1)=\emptyset.$

Now assume the relation $(4.4)$ holds true for $k$. Then $
L_{ij}^{[k]} = P_{elem}(v_{i}, v_{j}, k),$ that is $ L_{ij}^{[k]}$
represent the set of all elementary paths of length $k$ from
$v_{i}$ to $v_{j}$.

Using $(4.2)$ the $(i,j)$ entry of the matrix $L^{[k+1]} =
L\circ_{\ell} L^{[k]}$ is given explicitly by
\begin{equation}
L_{ij}^{[k+1]}=(L_{i1}\circ_{\ell} L_{1j}^{[k]})\cup \ldots \cup
(L_{im}\circ_{\ell} L_{mj}^{[k]})\cup\ldots \cup
(L_{in}\circ_{\ell} L_{nj}^{[k]}), ~~~k\geq 1. \label{(4.6)}
\end{equation}

We first evaluate the term $L_{im}\circ_{\ell} L_{mj}^{[k]}$ from
the equality $(4.6)$, for an fixed integer $m$ with $1\leq m\leq
n$. We have the following situations.

$\bullet~~~$ If $L_{im}=\varepsilon$ or $
L_{mj}^{[k]}=\varepsilon$, then $L_{im}\circ_{\ell}
L_{mj}^{[k]}=\varepsilon$, that is the set of elementary paths of
length $k+1$ from $v_{i}$ to $v_{j}$ stopping at vertex $v_{m}$ is
empty set.

$\bullet~~~L_{im}\neq\varepsilon$ and $
L_{mj}^{[k]}\neq\varepsilon$. By the induction hypothesis we
have\\[0.1cm]
$L_{mj}^{[k]}= P_{elem}(v_{m}, v_{j}, k)= \{ p_{mj}^{1},
p_{mj}^{2},\ldots, p_{mj}^{r-1},p_{mj}^{r} \},$ where $p_{mj}^{s}$
is an elementary path of length $k$ from $v_{m}$ to $ v_{j}$ for
$1\leq s\leq r$. Since $L_{im}=v_{i}v_{m}$ it follows
that\\[-0.1cm]
\[
L_{im}\circ_{\ell} L_{mj}^{[k]}=v_{i}v_{m}\circ_{\ell}\{
p_{mj}^{1}, p_{mj}^{2},\ldots, p_{mj}^{r-1},p_{mj}^{r} \}
\]
where $p_{mj}^{s}$ is regarded as a distinguished word, having the
set of vertices\\[0.1cm]
 $X_{mj}^{s}= \{v_{m}= v_{m_{0},j},
v_{m_{0}+1,j},\ldots, v_{m_{0}+k-1,j},v_{m_{0}+k,j}=v_{j} \}.$ The
element $ v_{i}v_{m}\circ_{\ell} p_{mj}^{s}~$ can be take the
following values:

$-~~~v_{i}v_{m}\circ_{\ell}
p_{mj}^{s}=\{v_{i}v_{m}v_{m_{0}+1,j}\ldots v_{m_{0}+k-1,j}v_{j}
\},$ if
$\{v_{i}\}\cap X_{mj}^{s}=\emptyset$,\\[0.1cm]
that is $v_{i}v_{m}\circ_{\ell} p_{mj}^{s}$ is an elementary path
of length $k+1$ from $v_{i}$ to $v_{j}$;

$-~~~v_{i}v_{m}\circ_{\ell} p_{mj}^{s}=\{\lambda\}, $ if
$\{v_{i}\}\cap X_{mj}^{s}\neq\emptyset$.

Then $L_{im}\circ_{\ell} L_{mj}^{[k]}=\{\lambda\} $ or
$~L_{im}\circ_{\ell} L_{mj}^{[k]} $ is a set which contains
$t~(1\leq t\leq s)$ elementary paths of length $k+1$ from $v_{i}$
to $v_{j}$ stopping at the vertex $v_{m}$.

Therefore, $ L_{ij}^{[k+1]} = \bigcup\limits_{m-1}^{n}
(L_{im}\circ_{\ell} L_{mj}^{[k]}) = P_{elem}(v_{i}, v_{j}, k+1) $.
Hence, $(4.4)$ holds for $k+1.$ This completes the inductive step
and proves the assertion in the case $i\neq j$.

{\bf (2)~Case $i=j$.}  If $L_{ii}^{[1]}= v_{i}v_{i},$  then $
C_{elem}(v_{i}, 1)= \{ (v_{i}, v_{i})\}$. Also, if
$L_{ii}=\varepsilon,$ then $C_{elem}(v_{i}, 1)=\emptyset.$ Hence,
$(4.5)$ holds for $k=1.$ Assume that $(4.5)$ holds true for $k$.
Then $ L_{ii}^{[k]} = C_{elem}(v_{i}, k),$ that is $ L_{ii}^{[k]}$
represent the set of all elementary circuits of length $k$
starting at $v_{i}$.

The $(i,i)$ entry of the matrix $ L^{[k+1]}$  is
$~L_{ii}^{[k+1]}=\bigcup\limits_{m=1}^{n}(L_{im}\circ_{\ell}
L_{mi}^{[k]}),~k\geq 1.$

We evaluate the term $L_{im}\circ_{\ell} L_{mi}^{[k]}$ for an
fixed integer $m$ with $1\leq m\leq n$. We have the following
situations.

$\bullet~~$ If $L_{im}=\varepsilon$ or $
L_{mi}^{[k]}=\varepsilon$, then $L_{im}\circ_{\ell}
L_{mi}^{[k]}=\varepsilon$.

$\bullet~~ L_{im}\neq\varepsilon$ and $
L_{mi}^{[k]}\neq\varepsilon$. Applying $(4.4)$ we have $~
L_{mi}^{[k]}= P_{elem}(v_{m}, v_{i}, k)= \{ q_{mi}^{1},
q_{mi}^{2},\ldots, q_{mi}^{r-1},q_{mi}^{r} \},$ where $q_{mi}^{s}$
is an elementary path of length $k$ from $v_{m}$ to $ v_{i}$ for
$1\leq s\leq r$. Since $L_{im}=v_{i}v_{m}$ it follows that
\[
L_{im}\circ_{\ell} L_{mi}^{[k]}=v_{i}v_{m}\circ_{\ell}\{
q_{mi}^{1}, q_{mi}^{2},\ldots, q_{mi}^{r-1},q_{mi}^{r} \}
\]
where $q_{mi}^{s}$ is regarded as a distinguished word, having the
set of vertices\\[0.1cm]
$Y_{mi}^{s}= \{v_{m}= v_{m_{1},i}, v_{m_{1}+1,i},\ldots,
v_{m_{1}+k-1,i},v_{m_{1}+k,i}=v_{i}\}$. The element $
v_{i}v_{m}\circ_{\ell} q_{mi}^{s} $ can be take the following
values:

$-~~~v_{i}v_{m}\circ_{\ell}
q_{mi}^{s}=\{v_{i}v_{m}v_{m_{1}+1,i}\ldots v_{m_{1}+k-1,i}v_{i}
\},$ if
$\{v_{i}\}\cap (Y_{mi}^{s}\setminus \{v_{i}\})=\emptyset$,\\[0.2cm]
that is $v_{i}v_{m}\circ_{\ell} q_{mi}^{s}$ is an elementary
circuit from $v_{i}$ to $v_{i}$ of length $k+1$;

$-~~~v_{i}v_{m}\circ_{\ell} q_{mi}^{s}=\{\lambda\}, $ if
$\{v_{i}\}\cap (Y_{mi}^{s}\setminus \{v_{i}\})\neq\emptyset$.

Then $L_{im}\circ_{\ell} L_{mi}^{[k]}=\{\lambda\} $ or
$~L_{im}\circ_{\ell} L_{mi}^{[k]} $ is a set which contains
$t_{1}~(1\leq t_{1}\leq s)$ elementary circuits of length $k+1$
from $v_{i}$ to $v_{i}$ stopping at $v_{m}$. Therefore, $
L_{ii}^{[k+1]} = C_{elem}(v_{i}, k+1) $. Hence, $(4.5)$ holds for
$k+1.$ This completes the inductive step and proves the assertion
in the case $i=j$. \hfill$\Box$

Applying Theorem 4.1 we give an answer of $({\bf EPP})$ and $({\bf
ECP})$  for a directed graph. In this purpose we give a new method
based on latin composition of distinguished languages, We will
called the {\it algorithm of latin composition of distinguished
languages} (shortly, {\bf LCDL}-algorithm).

{\it For a directed graph $G=(V, E)$ with $V=\{v_{1}, v_{2},
\ldots, v_{n}\}$, the {\bf LCDL}-algorithm consists from the
following steps:

{\bf Step 1.} Associate the latin matrix $L\in M(n, {\cal
P}^{\ast}(\Sigma_{dw}^{\ast}))$ to graph $G$;

{\bf Step 2.} For each $k (1\leq k\leq n)$ compute the matrix
$L^{[k]}=(L_{ij}^{[k]})$;

{\bf Step 3.}$~(i)$ For each pair $(i,j)$ and each $k (1\leq k\leq
n-1)$ enumerate the elementary paths of length $k$ from $v_{i}$ to
$v_{j}$ in $G$ {\rm (we apply the relation $(4.4)$)};

$(ii)~$For each $i$ and each $k (1\leq k\leq n)$ enumerate the
elementary circuits of length $k$ starting at $v_{i}$ in $G$ {\rm
(we apply the relation $(4.5)$).}}

\begin{Rem}
{\rm Let $G=(V,E)$ be a directed graph with $V=\{v_{1}, v_{2},
\ldots, v_{n}\}$.

$(i)$ {\bf (1)~} For $1\leq k\leq n-2$ and $i,j=\overline{1,n}$
with $i\neq j$, the elements $L_{ij}^{[k]}$ indicates the
elementary paths which are formed by $k$ arcs, so that
$L_{ij}^{[n-1]}$  determines all Hamiltonian paths in $G$ between
$v_{i}$ and $v_{j}$.

{\bf (2)~} For $1\leq k\leq n-1$ and $i=\overline{1,n}$, the
elements  $L_{ii}^{[k]}$ indicates the elementary circuits of
length $k$, so that $L_{ii}^{[n]}$ determines all Hamiltonian
circuits in $G$ starting at $v_{i}$.

$(ii)~$ The necessity to determine the elementary paths and
elementary circuits of maximum length in a directed graph arises.
 The {\bf LCDL}-algorithm determines the Hamiltonian paths (resp.,
Hamiltonian circuits) when there exist or determines the
elementary paths (resp., circuits) of maximum length when we have
no Hamiltonian paths (resp., circuits).} \hfill$\Box$
\end{Rem}
\begin{Rem}
(\cite{bagu}) {\rm The powers of the adjacency matrix $A\in
M(n;{\bf N})$ associated to graph $G=(V,E)$ are used to find the
number of distinct paths and distinct circuits between two
vertices in $V$ (two paths or circuits of length $k$ are distinct
if they visit a different sequence of vertices). More precisely:}

Let $A^{k}=(A_{ij}^{k})$ be the $k-$th power of the adjacency
matrix $A$. Then:\\[-0.1cm]
\[
A_{ij}^{k}= |P(v_{i}, v_{j},k)| ~~~\hbox{and}~~~ A_{ii}^{k}=
|C(v_{i}, k)|~~~\hbox{for all}~~k\geq 1.~~~
 \hfill\Box
 \]
\end{Rem}
\begin{Ex}
{\rm Let $G=(V, E)$ be a directed graph with vertex set
$V=\{v_{1}, v_{2}, v_{3}, v_{4}\}$ and the adjacency matrix $A\in
M_{4}({\bf N})$ of $G$ where\\[-0.2cm]
\[
A = \left( \begin{array}{cccc}
 1 & 1 & 1 & 1 \\
 0 & 1 & 1 & 1 \\
 0 & 0 & 0 & 1 \\
 0 & 0 & 0 & 0
\end{array} \right).
\]
\\[-0.2cm]
The geometric representation of graph $G$ is given in Figure
4.1.\\[-0.5cm]
\begin{center}
\includegraphics[width=7cm]{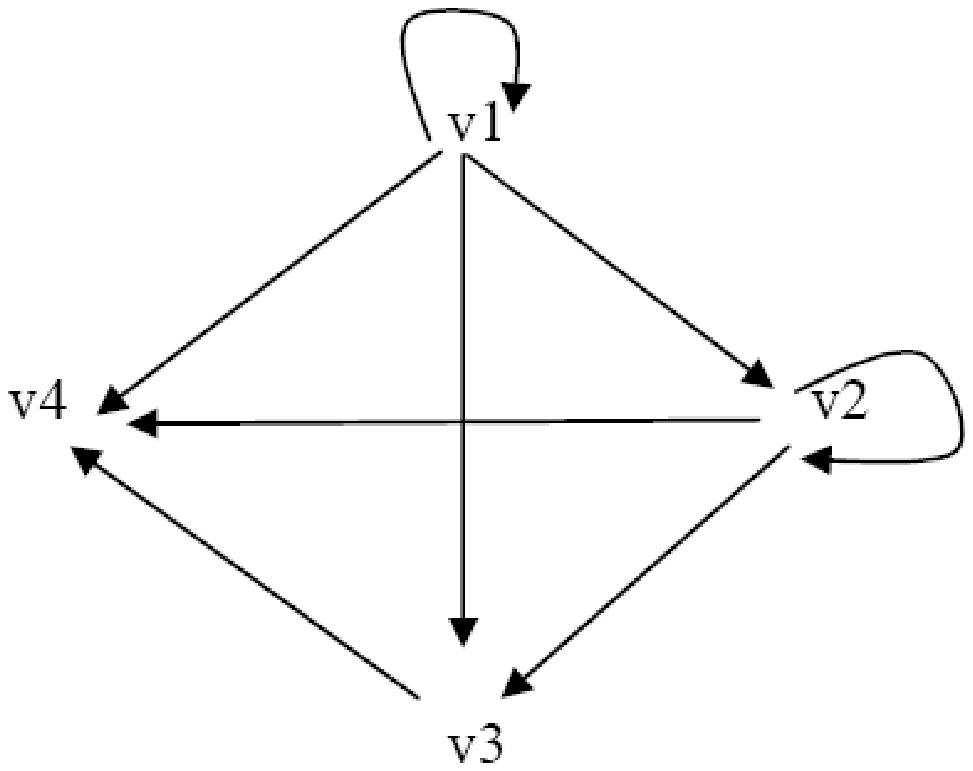}
\end{center}
\hspace*{6.5cm} Fig.4.1.\\

Computing the power $A^{3}$ of the adjacency matrix $A$ in terms
of numbering semiring $({\bf N}, +, \cdot, 0, 1)$ we have
$A_{14}^{3}=5~$ and $A_{22}^{3}=1$. Then:

${\bullet~~}$ Between $v_{1}$ and $v_{4}$ there exist $5$ paths of
length $3$, namely: $p_{1}^{3}=( v_{1}, v_{2},
v_{3},v_{4}),$\\[0.1cm]
$p_{2}^{3}=( v_{1}, v_{1}, v_{1},  v_{4}),~ p_{3}^{3}=(v_{1},
v_{1}, v_{2}, v_{4}),~p_{4}^{3}= (v_{1}, v_{2}, v_{2}, v_{4}),~
p_{5}^{3}=( v_{1}, v_{1}, v_{3}, v_{4})$;

${\bullet~~} G$ has one circuit of length $3$ starting at $v_{2}$,
namely  $ c_{1}^{3}=( v_{2}, v_{2}, v_{2}, v_{2})$.

$(ii)~$ We consider the alphabet $\Sigma = \{ v_{1}, v_{2}, v_{3},
v_{4} \} $. For to find the elementary paths and elementary
circuits  in graph $G$ we apply {\bf LCDL}-algorithm.

The latin matrix $L\in M_{4}({\cal P}(\Sigma_{dw}^{\ast}))$
associated to graph $G$ is\\[-0.2cm]
\[
L = \left( \begin{array}{cccc}
 v_1v_1 & v_1v_2 & v_1v_3& v_1v_4 \\
 \varepsilon & v_2v_2 & v_2v_3 & v_2v_4 \\
 \varepsilon & \varepsilon & \varepsilon & v_3v_4 \\
 \varepsilon & \varepsilon & \varepsilon &  \varepsilon
\end{array} \right).
\]
\\[-0.2cm]
We compute the powers od the latin matrix $L$ in terms of the
semiring of distinguished languages. We first compute the matrix
$L^{[2]}=(L_{ij}^{[2]}).$ We have
\[
L^{[2]}=L\circ_{\ell} L= =\left( \begin{array}{cccc}
 \varepsilon & \varepsilon & v_1v_2v_3& \{v_{1} v_{2}v_{4}, v_{1} v_{3}v_{4} \} \\
 \varepsilon & \varepsilon & \varepsilon & v_{2}v_{3}v_{4} \\
 \varepsilon & \varepsilon & \varepsilon & \varepsilon \\
 \varepsilon & \varepsilon & \varepsilon &  \varepsilon
\end{array} \right).
\]

For example,  $L_{14}^{[2]}$ is computed as
follows\\[0.1cm]
$L_{14}^{[2]}= (L_{11}\circ_{\ell}L_{14})\cup
(L_{12}\circ_{\ell}L_{24})\cup(L_{13}\circ_{\ell}L_{34})\cup(L_{14}\circ_{\ell}L_{44})=(v_{1}v_{1}\circ_{\ell}v_{1}v_{4})\cup$\\[0.1cm]
$\cup(v_{1}v_{2}\circ_{\ell}v_{2}v_{4})\cup(v_{1}v_{3}\circ_{\ell}v_{3}v_{4})\cup
(v_{1}v_{4}\circ_{\ell}\varepsilon) = \varepsilon\cup
\{v_{1}v_{2}v_{4}\}\cup \{v_{1}v_{3}v_{4}\}\cup
\varepsilon=\{v_{1}v_{2}v_{4}, v_{1}v_{3}v_{4}\}.$\\[-0.2cm]

Since $L_{14}^{[2]}= \{v_{1}v_{2}v_{4}, v_{1}v_{3}v_{4}\}$ it
follows that there exist two elementary paths of length $2$ from
$v_{1}$ to $v_{4}$ and we have$~ P_{elem}(v_{1}, v_{4}, 2)=\{
(v_{1}, v_{2}, v_{4}), ~(v_{1}, v_{3}, v_{4})\}.$

The matrices $ L^{[3]}=L\circ_{\ell} L^{[2]}$ and $
L^{[4]}=L\circ_{\ell} L^{[3]}$ are given by
\[
L^{[3]}=\left( \begin{array}{cccc}
 \varepsilon & \varepsilon & \varepsilon & v_{1} v_{2}v_{3}v_{4} \\
 \varepsilon & \varepsilon & \varepsilon & \varepsilon \\
 \varepsilon & \varepsilon & \varepsilon & \varepsilon \\
 \varepsilon & \varepsilon & \varepsilon &  \varepsilon
\end{array} \right)~~~\hbox{and}~~~L^{[4]}=\left( \begin{array}{cccc}
 \varepsilon & \varepsilon & \varepsilon & \varepsilon \\
 \varepsilon & \varepsilon & \varepsilon & \varepsilon \\
 \varepsilon & \varepsilon & \varepsilon & \varepsilon \\
 \varepsilon & \varepsilon & \varepsilon &  \varepsilon
\end{array} \right).
 \]
Using the matrices $L^{[3]}$ and $L^{[4]}$ one obtains the
following results:

$\bullet~~ L_{14}^{[3]}=\{v_{1}v_{2}v_{3}v_{4}\}$. Then
$P_{elem}(v_{1}, v_{4}, 3)=\{(v_{1}, v_{2}, v_{3}, v_{4})\}$ and
$G$ has only one Hamiltonian path. The set of elementary paths of
maximum length from $v_{2}$ to $v_{4}$ is  $P_{elem}(v_{2}, v_{4},
2)=\{ ( v_{2}, v_{3}, v_{4})\}$, since
$L_{24}^{[2]}=\{v_{2}v_{3}v_{4}\}$ and $L_{24}^{[3]}=\varepsilon.$

$\bullet~~ L_{ii}^{[k]}=\varepsilon~$ for $2\leq k\leq 4$ and
$i=\overline{1,4}.$ Then $C_{elem}(v_{i}, v_{i}, k)=\emptyset $.
 The elementary circuits of maximum length are those of length
$1$.}\hfill$\Box$
\end{Ex}

{\bf Application.  The finding of Hamiltonian paths and
Hamiltonian circuits of  minimal cost in a weighted directed
graph}.

{\it The {\bf LCDL}-algorithm can be to use for solving of the
following two problems in a weighted directed graph $G=(V,E)$:

$(i)~$ find a Hamiltonian path of minimal cost between two
vertices in $G$;

$(ii)~$ for each  $v\in G$, find a Hamiltonian circuit of minimal
cost starting at $v$.}

One way to solve the above problems consists of searching all
possible Hamiltonian paths and Hamiltonian circuits (we apply {\bf
LCDL}-algorithm) and computing their cost (we apply the relation
(3.1)).
\begin{Ex}
{\rm Let $G=(V, E, w)$ be a weighted directed graph where $V=\{1,
2, 3, 4, 5\},$\\[0.1cm]
$E=\{ (1, 2), (1, 3), (1, 5), (2, 1), (2, 5), (3, 2), (4, 3), (4,
5), (5, 1), (5, 2), (5, 3), (5, 4) \}$ and the cost
function $ w_{cost}: E\to {\bf R},~(i,j)\rightarrowtail w_{cost}(i,j)=w_{ij}$ given by\\[0.2cm]
$\begin{array}{|c|c|c|c|c|c|c|c|c|c|c|} \hline (i,j)& (1, 2)&
(1,3)&(1, 5)&(2, 1)&(2, 5)&(3, 2)&(4, 3)&(4, 5)&(5, 1)&(5, 2)\\
\hline
 w_{ij}& 4 & 2 & 6 & 3 & 3 & 1 & 5 & 4 & 6 & 1 \cr \hline
\end{array}$\\
$\begin{array}{|c|c|c|} \hline (i,j)& (5, 3)&(5, 4)\\
\hline
 w_{ij}& 2 &1 \cr \hline
\end{array}$\\

 The latin matrix $L\in M_{5}({\cal P}(\Sigma_{dw}^{\ast}))$
($\Sigma = \{ 1, 2, 3, 4, 5 \} $) associated to graph $G$ is
\[
L = \left( \begin{array}{ccccc}
\varepsilon & 12 & 13 & \varepsilon & 15\\
 21 & \varepsilon & \varepsilon &\varepsilon & 25  \\
 \varepsilon & 32 & \varepsilon & \varepsilon & \varepsilon \\
 \varepsilon& \varepsilon & 43 & \varepsilon & 45\\
 51 & 52 & 53 & 54 & \varepsilon\\
 \end{array} \right).
\]

Let us we compute the powers of the latin matrix $L$. We have
\[
L^{[2]}=L\circ_{\ell} L =\left( \begin{array}{ccccc}
 \{121, 151\}& \{132,152\}& 153& 154 & 125 \\
 251 & \{ 212, 252\}& \{213, 253\} & 254 & 215 \\
321 & \varepsilon & \varepsilon & \varepsilon & 325 \\
 451 & \{ 432, 452\}& 453 &  454 & \varepsilon\\
521 & \{ 512, 532\}& \{513, 543\} & \varepsilon &\{515, 525, 545\} \\
\end{array} \right).
\]

The matrix $ L^{[3]} = L\circ_{\ell} L^{[2]}$ has the following
lines:

$L_{1}^{[3]}:~~(\{1321,1521,1251\},~1532,~ \{1543,1253\},~1254,~1325);$\\[-0.2cm]

$L_{2}^{[3]}:~~(\varepsilon,~ \{2512,2132,2532,2152\},~  \{
2513,2543,2153\},~  2154, ~~
 \varepsilon );$ \\[-0.2cm]

$L_{3}^{[3]}:~~(3251,~~ \varepsilon,~ \{3213,3253\},~3254,~ 3215  );$ \\[-0.2cm]

$L_{4}^{[3]}:~~(\{4321,4521\},~\{4512, 4532\},~ 4513,~~ \varepsilon,~~ 4325 );$ \\[-0.2cm]

$L_{5}^{[3]}:~~(5321,~ \{5132,5432\},~ 5213,~ \varepsilon,~ \{5215,5125,5325\}).$\\[-0.2cm]

The matrix $ L^{[4]} = L\circ_{\ell} L^{[3]}$ has the following
lines:

$L_{1}^{[4]}:~~(\{15321,13251\},~  15432, ~12543,~  13254,~~  \varepsilon );$\\[-0.2cm]

$L_{2}^{[4]}:~~(\varepsilon,~ \{25132,25432,21532\},~ 21543,~~  \varepsilon,~~  \varepsilon);$\\[-0.2cm]

$L_{3}^{[4]}:~~(\varepsilon,~~ \varepsilon,~ \{32513,32543,32153\},~~ \varepsilon,~~  \varepsilon);$\\[-0.2cm]

$L_{4}^{[4]}:~~(\{45321,43251\},~ 45132,~  45213,~  43254,~  43215);$\\[-0.2cm]

$L_{5}^{[4]}:~~(54321,~~ \varepsilon,~~  \varepsilon,~~ \varepsilon,~ \{53215,51325,54325\}).$\\[0.1cm]

Finally, we have
\[
L^{(5)} = L^{(4)} \circ_{\ell} L= \left( \begin{array}{ccccc}
154321 &  \varepsilon & \varepsilon & \varepsilon & \varepsilon\\
\varepsilon & 215432 & \varepsilon & \varepsilon & \varepsilon\\
\varepsilon & \varepsilon & \varepsilon & \varepsilon & \varepsilon\\
\varepsilon & \varepsilon & \varepsilon & 432154 & \varepsilon\\
\varepsilon & \varepsilon & \varepsilon &  \varepsilon & 543215
\end{array} \right).
\]

$\bullet ~~$ From $L^{[4]}$ it follows that $G$ has $11$ Hamiltonian
paths. For example,
 $~P_{elem}(4,1,4)=\{p_{1,H}^{4}=(4,5,3,2,1), p_{2,H}^{4}=(4,3, 2,
5,1)\}$ since $L_{41}^{[4]}=\{45321, 43251\}$. We have
$~w(p_{1,H}^{4})= 10$ and $~w(p_{2,H}^{4})=15.~$ It follows that
$p_{2,H}^{4}$ is a Hamiltonian path between the vertices $4$ and $1$
having the maximal cost equal to $15$.

$\bullet ~~$ From $L^{[5]}$ it follows that  $G$ has $4$ Hamiltonian
circuits.  For example,
$C_{elem}(1,1,5)=L_{11}^{[5]}=\{c_{1,H}^{5}=(1, 5, 4, 3, 2, 1)\}$.
We have $~w(c_{1,H}^{5}) = 16$. Hence $c_{1,H}^{5}$ is a Hamiltonian
circuit starting at vertex $1$ having the maximal  cost equal to
$16$.}\hfill$\Box$
\end{Ex}
{\bf Conclusions}. The {\bf LCDL}- algorithm can be easily
programmed and gives an efficient solution of {\bf EPP} and {\bf
ECP} for a finite directed graph. This can be seen as an improved
version of  Kaufmann's algorithm.

Let us list some classes of algebraic path problems, which can be
reduced to applying of the {\bf LCDL}-algorithm:

$(i)~$ enumeration of the elementary paths (resp., circuits);

$(ii)$ determination of the elementary paths (resp., circuits) of
maximum length;

$(iii)$ testing a graph $G$ for having Hamiltonian paths or
Hamiltonian circuits;

$(iv)$ optimization (Hamiltonian path or Hamiltonian circuit of
minimal cost). \hfill$\Box$

Author's adress\\

West University of Timi\c soara,\\
Department of Mathematics, Bd. V. P{\^a}rvan, no. 4, 300223, Timi\c soara, Romania\\
\hspace*{0.7cm} E-mail: ivan@math.uvt.ro\\
\end{document}